\documentclass{amsart}
\usepackage{amssymb}
\usepackage{amsmath}
\usepackage{amsfonts}
\usepackage{pstricks}
\usepackage{pst-node}
\usepackage{graphicx}
\theoremstyle{plain}
\newtheorem{thm}{Theorem}[section]

\newtheorem{dfn}[thm]{Definition}

\newtheorem{qws}[thm]{Question}

\theoremstyle{remark} 

\def\pmc#1{\setbox0=\hbox{#1}
    \kern-.1em\copy0\kern-\wd0
    \kern.1em\copy0\kern-\wd0}

\begin{document}

\title[A new class of homology and cohomology 3-manifolds]{A new class of homology and \\ cohomology 3-manifolds}
\author{D.J. Garity, U.H. Karimov, D.D. Repov\v{s}, and F. Spaggiari}

\address{
Department of Mathematics,
Oregon State University,
Corvallis, Oregon,
USA}
\email{garity@math.oregonstate.edu}

\address{
Institute of Mathematics, 
Academy of Sciences of Tajikistan,  
Dushanbe, 
Taji\-ki\-stan}
\email{umedkarimov@gmail.com}

\address{
Faculty of Education and 
Faculty of Mathematics and Physics,
University of Ljubljana
\&
Institute of Mathematics, Physics and Mechanics, 
Ljubljana, 
Slovenia}
\email{dusan.repovs@guest.arnes.si}

\address{
Department of Mathematics, 
University of Modena and Reggio Emilia, 
Mo\-de\-na,
Italy}
\email{fulvia.spaggiari@unimore.it}


\begin{abstract}
We show that for any set of primes 
$\mathcal{P}$
there exists a space 
$M_{\mathcal{P}}$ which
is 
a homology and cohomology 3-manifold with coefficients in $\mathbb{Z}_{p}$ for
$p\in \mathcal{ P}$  
and
is not  a homology 
or cohomology 3-manifold  with coefficients in $\mathbb{Z}_q$  for $q\not\in \mathcal{ P}$.
Moreover,
$M_{\mathcal{P}}$ is neither a  homology nor cohomology 3-manifold with coefficients in 
$\mathbb{Z}$ or  $\mathbb{Q}$.
\end{abstract}

\maketitle

\section{Introduction}

In 1908 Tietze \cite{T} constructed his famous $3$-manifolds
$L(p, q)$ called {\it lens spaces}. These spaces have many
interesting properties. For example,
lens spaces 
$L(5,1)$ and
$L(5,2)$ have isomorphic fundamental groups and the same homology groups,
but they do not have the same homotopy type (this was first proved in 1919 by Alexander \cite{A}). 

It is also well known that
for every prime $q,$
the lens space $M^3 = L(q, 1)$ has the following properties:
\begin{itemize}
\item
$M^3$
is a $\mathbb{Z}_p$-homology 3-sphere 
for every prime $p \not= q$;
\item $M^3$
is not a $\mathbb{Z}_q$-homology 3-sphere; and
\item $M^3$
is not a $\mathbb{Z}$-homology 3-sphere.
\end{itemize}
We shall generalize
this classical
result as follows:

\begin{thm}\label{Thm:Main}
Given any set of primes $\mathcal{P},$ there exists a 3-dimensional compact metric space $M_{\mathcal{P}}$ which
has the following properties:
\begin{itemize}
\item
$M_{\mathcal{P}}$
is a 
${\mathbb{Z}_p}$-homology 
3-manifold 
  for every
$p\in \mathcal{ P}$;
\item
$M_{\mathcal{P}}$
is a 
${\mathbb{Z}_p}$-cohomology 
3-manifold 
  for every
$p\in \mathcal{ P}$;
\item
$M_{\mathcal{P}}$
 is not a 
${\mathbb{Z}_q}$-homology 3-manifold  for any $q\not\in
\mathcal{ P}$;
\item
$M_{\mathcal{P}}$
 is not a 
${\mathbb{Z}_q}$-cohomology 3-manifold  for any $q\not\in
\mathcal{ P}$;
\item
 $M_{\mathcal{P}}$ is not  a 
 $\mathbb{Z}$-homology 3-manifold; 
\item
 $M_{\mathcal{P}}$ is not  a 
 $\mathbb{Z}$-cohomology 3-manifold; 
\item 
$M_{\mathcal{P}}$
is not a $\mathbb{Q}$-homology 3-manifold; and
\item 
$M_{\mathcal{P}}$
is not a $\mathbb{Q}$-cohomology 3-manifold.
\end{itemize}
\end{thm}


\section{Preliminaries}

First, we fix the terminology, notation, and recall some well
known facts.   We let $L$ be
the ring of integers $\mathbb{Z}$ or a
field.

\begin{dfn}\label{chm} 
{\rm (cf. \cite[Corollary V.16.9]{Br})} 
A   locally compact
space $X$ is called a cohomology manifold of dimension $n$ 
over $L$ {\rm (}denoted as $n-cm_L${\rm )} if
\begin{enumerate}
  \item $X$ has finite
cohomological dimension over $L$
$\rm{ ( dim}_L X<\infty \rm{)}$;
  \item  $X$ is  cohomologically
locally connected over $L$ {\rm (}$clc^{\infty}_L${\rm )}; and
  \item for all $x \in X,$  
$$\check{H}^p(X, X- \{x\} ; L) \cong
\left\{
\begin{array}{ll}
L     & for \ \  p=n \\
0 & for \ \ p \neq n\\
\end{array}
\right. $$ 
 where $\check{H}^*$ are \v{C}ech
cohomology groups with coefficients in  $L.$
\end{enumerate}
 \end{dfn}

\begin{dfn} A locally compact space $X$ is called a homology manifold of
dimension $n$ over 
$L$ {\rm (}denoted as $n-hm_L${\rm )} if
\begin{enumerate}
\item $X$  has finite cohomological
dimension over $L$
$\rm{ ( dim}_L X<\infty \rm{)}$; 
 \item for all $x \in X,$  
$${H}_p(X, X- \{x\} ; L) \cong
\left\{
\begin{array}{ll}
L     & for \ \  p=n \\
0 & for \ \ p \neq n\\
\end{array}
\right. $$ 
 where ${H}_*$ are Borel-Moore homology groups with coefficients in  $L.$
\end{enumerate}
Homology
manifolds will stand for homology $\mathbb{Z}-$manifolds.
\end{dfn}

Any $n-$dimensional cohomology manifold $(n-cm_L)$ is
also
an $n-$dimensional homology manifold $(n-hm_L)$ by \cite[Theorem V.16.8]{Br}.
Therefore we will construct only cohomology manifolds which will
be automatically homology manifolds by this theorem.

{\bf Example 2.3.}
For the construction and some simple properties of lens spaces see
\cite{H, ST}. In particular, the homology groups
of the lens space $M^3 =
L(q, 1)$ are

$$H_n(M^3; \mathbb{Z}) =
\left\{
\begin{array}{ll}
\mathbb{Z}     & n=0,3 \\
\mathbb{Z}_{q} & n=1   \\
0       & n=2 \text{ or }  n\ge 4
\end{array}
\right. $$

By the Universal Coefficients Theorem, we have for any abelian group $G,$  

$$H_n(M^3; G)\cong H_n(M^3; \mathbb{Z})\otimes G \oplus H_{n-1}(M^3; \mathbb{Z}) \ast G.$$

\ \\
Therefore, if $p$ and $q$ are prime and $p\neq q$ then 

$$H_n(M^3; \mathbb{Z}_p) \cong
\left\{
\begin{array}{ll}
\mathbb{Z}_p      & \text{if }n=0,3 \\
0          & \hbox{otherwise}
\end{array}
\right.
$$

whereas

$$H_n(M^3; \mathbb{Z}_q) \cong
\left\{
\begin{array}{ll}
\mathbb{Z}_{q} & n = 0, 1, 3   \\
0       & \hbox{otherwise}
\end{array}
\right. $$

Consider now the suspension
$N^4=\Sigma M^3$
of $M^3$. Local conditions of 
Definitions 2.1 and 2.2 are clearly satisfied for $N^4$
since $M^3$ is a 
3-manifold. Therefore $N^4$
 is a
4-$hm_{\mathbb{Z}_p}$ and
a
4-$cm_{\mathbb{Z}_p}$.
However, $N^4$
is neither
 a 4-$hm_{\mathbb{Z}_q}$
nor
a
4-$cm_{\mathbb{Z}_q}$ if $p$ and $q$ are prime and $p
\not= q$ (cf. \cite{Br}). 

\section{Proof of Theorem \ref{Thm:Main}}

Let 
$\mathcal{P} = \{p_i\}_{i \in K}$, 
for 
$K = {\mathbb N}$ 
or 
$K = \{1,\dots, k\}$, 
be a set of some prime
numbers. 
If the set $K$ is infinite then we define the numbers
$n_i$ as $n_i = p_1\cdot p_2\cdot p_3\cdots p_i.$
If the set $K$ is finite and consists of exactly $k$ elements, then
define $n_i$ as $n_i = p_1\cdot p_2\cdot p_3\cdots p_k$ for all
$i.$

Let $X$ be a solenoid in the 3-dimensional sphere $S^3$, i.e. the inverse limit of solid tori corresponding to the
following inverse system:
$$\mathbb{Z} \stackrel{n_1}{\longleftarrow} \mathbb{Z} \stackrel{n_2}
{\longleftarrow} \mathbb{Z} \stackrel{n_3}{\longleftarrow}\ldots
$$
naturally embedded in $S^3,$ see e.g. \cite[Chapter IX, Exercise
E.4]{ES}.

Let us prove that the quotient space $S^3/X$ is a cohomology
3-manifold $cm_{\mathbb{Z}_p}$. It is obvious that $S^3/X$ is
$3-$dimensional, compact and metrizable. So the space $S^3/X$ satisfies
the condition (1) of Definition \ref{chm}.

To prove that the space $S^3/X$ satisfies
the conditions (2) and (3) of Definition \ref{chm}, let's calculate 
first the groups $\check{H}^n(S^3/X, \{x\};
G)$ with respect to
the
one-point subset $\{x\}=X/X$  for $G \cong \mathbb{Z} _p, p\in \mathcal{P};\,\, G \cong \mathbb{Z} _q, q\not\in  \mathcal{P};\,\, G \cong \mathbb{Z};\,\, G \cong \mathbb{Q}.$ Since $S^3/X$ is connected and $3-$dimensional it follows
that
\begin{equation}\label{X^{ast3}}
\check{H}^0(S^3/X, \{x\}; G) \cong 0 \ \ \hbox{ and} \ \ 
\check{H}^n(S^3/X, \{x\}; G) \cong 0 \ \ \hbox {for} \ n > 3.
\end{equation}

Let $U_i$ be the open $i$-th solid torus neighborhood of $X$ in $S^3$
(see \cite{ES}).
Then $\{U_i/X\}_{i\in \mathbb{N}} $ is a
neighborhood base  of $x$ in $S^3/X.$
By continuity of the {\v C}ech cohomology and by the Excision Axiom
it follows that:
$$\check{H}^n(S^3/X, \{x\}; G) \cong \check{H}^n(\lim_{\leftarrow}S^3/\overline{U}_i,
\overline{U}_i/\overline{U}_i; G)\cong$$ $$\cong
\lim_{\to}\check{H}^n(S^3/\overline{U}_i,
\overline{U}_i/\overline{U}_i; G)\cong
\lim_{\to}\check{H}^n(S^3,
 \overline{U}_i; G).$$

For $n = 1$ we have the exact sequence of the pair $(S^3, \overline{U_i})$:
$$\check{H}^0(S^3; G)\longrightarrow \check{H}^0(\overline{U_i};
G) \longrightarrow \check{H}^1(S^3, \overline{U}_i;
G)\longrightarrow \check{H}^1(S^3, G)\longrightarrow \check{H}^1(\overline{U_i};
G).$$
Since the $1-$dimensional cohomology of the $3-$sphere is trivial for any
group of coefficients $G$ and $\overline{U_i}$ is connected space for every $i,$ it follows that $\check{H}^1(S^3,
\overline{U}_i; G) \cong 0,$ therefore $\check{H}^n(S^3/X, \{x\}; G) \cong 0$ and, in particular, 
\begin{equation}\label{H^1}
\check{H}^1(S^3/X, \{x\}; \mathbb{Z}_p) \cong 0
\end{equation}
and
\begin{equation}\label{H^1ZQZ_g}
\check{H}^1(S^3/X, \{x\}; \mathbb{Z}_q) \cong 0,
\check{H}^1(S^3/X, \{x\}; \mathbb{Z}) \cong 0, \check{H}^1(S^3/X,
\{x\}; \mathbb{Q}) \cong 0.
\end{equation}

For $n = 2$ we have the following exact sequence of the pair $(S^3, \overline{U_i})$:
$$\check{H}^1(S^3; G)\longrightarrow \check{H}^1(\overline{U}_i;
G) \longrightarrow \check{H}^2(S^3, \overline{U}_i;
G)\longrightarrow \check{H}^2(S^3, G)\longrightarrow \check{H}^2(\overline{U_i};
G).$$

The cohomology groups $\check{H}^1(S^3; G)$ and $\check{H}^2(S^3, G)$ are
trivial, and the groups $\check{H}^1(\overline{U}_i; G)$ are
isomorphic to $G$ since $U_i$ has the homotopy type of a circle.
The homomorphisms $\check{H}^1(\overline{U}_i; G)\rightarrow
\check{H}^1(\overline{U}_{i+1}; G)$ are multiplications by $n_i$
that take the group $G$  into itself. 

Therefore for the group of coefficients  $G \cong \mathbb{Z}_p$ it
follows that
\begin{equation}\label{H^2}
\check{H}^2(S^3/X, \{x\}; \mathbb{Z}_p) \cong 0.
\end{equation}

\eject

However, 

\begin{equation}\label{H^2alt}
\check{H}^2(S^3/X, \{x\}; \mathbb{Z}_q)  \ncong 0 , 
\check{H}^2(S^3/X,
\{x\}; \mathbb{Z}) \ncong 0, 
\check{H}^2(S^3/X, \{x\}; \mathbb{Q})
\ncong 0.
\end{equation}

For $n = 3$ consider  the following cohomology exact sequence for the pair $(S^3, \overline{U_i})$:
$$\check{H}^2(\overline{U}_i;
G) \longrightarrow \check{H}^3(S^3, \overline{U}_i;
G)\longrightarrow \check{H}^3(S^3, G)\longrightarrow
\check{H}^3(\overline{U}_i; G).$$

Since $\overline{U}_i\simeq S^1$, it follows that:
\begin{equation}\label{H^3}
\check{H}^3(S^3/X, x; G) \cong G.
\end{equation}

Let us calculate the groups $\check{H}^n(S^3/X - \{x\}; \mathbb{Z}_p)$. The
space $S^3/X - \{x\}$ is the union $\bigcup_{i=1}^{\infty}(S^3 -
U_i)$ of an
increasing sequence of ``complementary'' solid tori.

For $n = 1$ we have the following exact sequence of Milnor-Harlap
\cite[Theorem 1]{Kh}:
$$0 \rightarrow \underleftarrow{\lim}^{(1)} \check{H}^{0}(S^3 - U_i; \mathbb{Z}_p)
\rightarrow \check{H}^1(S^3 - X; \mathbb{Z}_p)\rightarrow
\underleftarrow{\lim}\check{H}^1(S^3 - U_i; \mathbb{Z}_p)\rightarrow 0,$$
where $\underleftarrow{\lim}^{(1)}$ is the first derived functor of
the functor of the inverse limit. Since $p\in \mathcal{P}$ it
follows that the inverse limit $\underleftarrow{\lim}\check{H}^1(S^3 -
U_i; \mathbb{Z}_p)$ is trivial. The group $\underleftarrow{\lim}^{(1)}
\check{H}^{0}(S^3 - U_i; \mathbb{Z}_p)$ is trivial since the corresponding
inverse sequence satisfies the Mittag-Leffler (ML) condition, so we have

\begin{equation}\label{H^1(S3-x)}
\check{H}^1(S^3 - X; \mathbb{Z}_p) \cong 0.
\end{equation}

Analogously, it is easy to see that
\begin{equation}\label{H^1(S3-x)ZQZ_p}
\check{H}^1(S^3 - X; \mathbb{Z}_q) \ncong 0 \ \hbox
{for}\ \ q \notin \mathcal{P},\ \
\check{H}^1(S^3 - X; \mathbb{Z}) \cong 0,\ \check{H}^1(S^3 - X;
\mathbb{Q}) \ncong 0.
\end{equation}

Let $n = 2,$ then we have the Milnor-Harlap exact sequence for the presentation 
$S^3/X -  \{x\} = \bigcup_{i=1}^{\infty}(S^3 -
U_i)$:
$$0 \rightarrow \underleftarrow{\lim}^{(1)} \check{H}^{1}(S^3 - U_i; \mathbb{Z}_p)
\rightarrow \check{H}^2(S^3 - X; \mathbb{Z}_p)\rightarrow
\underleftarrow{\lim}\check{H}^2(S^3 - U_i; \mathbb{Z}_p)\rightarrow 0.$$
The groups $\underleftarrow{\lim}^{(1)} \check{H}^{1}(S^3 -
U_i; \mathbb{Z}_p)$ are trivial since the groups $\check{H}^{1}(S^3 -
U_i; \mathbb{Z}_p)$ are isomorphic to the finite group $\mathbb{Z}_p$ and the corresponding
inverse sequence satisfies the ML condition. The
groups $\check{H}^2(S^3 - U_i; \mathbb{Z}_p)$ are also trivial since the "complementary" solid tori have the homotopy type of the circle. Therefore
\begin{equation}\label{H^2(S3-x)}
\check{H}^2(S^3 - X; \mathbb{Z}_p) \cong 0.
\end{equation}

For $n = 3$ we have the exact sequence of Milnor-Harlap for the same presentation of $S^3/X  -  \{x\}$ as before:
$$0 \rightarrow \underleftarrow{\lim}^{(1)} \check{H}^{2}(S^3 - U_i; \mathbb{Z}_p)
\rightarrow \check{H}^3(S^3 - X; \mathbb{Z}_p)\rightarrow
\underleftarrow{\lim}\check{H}^3(S^3 - U_i; \mathbb{Z}_p)\rightarrow 0.$$
The groups $\check{H}^3(S^3 - U_i; \mathbb{Z}_p)$ and $\check{H}^{2}(S^3 -
U_i;\mathbb{Z}_p)$ are trivial since the spaces $S^3 - U_i$ have the homotopy
type of a circle. Therefore:
\begin{equation}\label{H^3(S3-x)}
\check{H}^3(S^3 - X; \mathbb{Z}_p) \cong 0.
\end{equation}

Next, let us calculate the groups $\check{H}^n(S^3/X, S^3/X -
X/X; G)$ for certain groups $G.$

Since the space $S^3/X$ is connected and $\dim S^3/X = 3$ it
follows that these groups are trivial groups for $n = 0, n > 3.$

Since the space $S^3 - X$ is connected and $\check{H}^1(S^3/X;
\mathbb{Z}_p) \cong 0$ by (\ref{H^1}),
it follows by the exact cohomology
sequence of the pair 
$(S^3/X, S^3/X - X/X)$ (or the pair $(S^3/X, S^3  -  X)$ since
$S^3/X  -  X/X \cong S^3  -  X$) that
\begin{equation}\label{cmZ_p1}
\check{H}^1(S^3/X, S^3 - X; \mathbb{Z}_p) \cong 0.
\end{equation}
By the exact sequence: $$\check{H}^{1}(S^3 - X; \mathbb{Z}_p)
\stackrel{\delta}\longrightarrow  \check{H}^2(S^3/X, S^3 - X; \mathbb{Z}_p)
\longrightarrow  \check{H}^2(S^3/X; \mathbb{Z}_p)\longrightarrow \check{H}^{2}(S^3 - X; \mathbb{Z}_p)
$$ and since the groups $\check{H}^{1}(S^3 - X; \mathbb{Z}_p)$ and $\check{H}^2(S^3/X;
\mathbb{Z}_p)$ are trivial by (\ref{H^1(S3-x)}) and (\ref{H^2}) it follows 
that
\begin{equation}\label{cmZ_p2}
\check{H}^2(S^3/X, S^3 - X; \mathbb{Z}_p) \cong 0.
\end{equation}
For any group of coefficients the corresponding homomorphism
$\delta$ is a monomorphism by (\ref{H^1ZQZ_g}). Since the
groups $\check{H}^{1}(S^3 - X; \mathbb{Z}_q)\
\hbox{for} \ q\notin \mathcal{P}, \ \check{H}^{1}(S^3 - X; Q) $ are nontrivial 
by (\ref{H^1(S3-x)ZQZ_p}),
and the groups $\check{H}^{1}(S^3/X; \mathbb{Z}_q), \check{H}^{1}(S^3/X; Q)$ are trivial,
it follows that
\begin{equation}\label{cmZ_QZ_p}
\check{H}^2(S^3/X, S^3 - X; \mathbb{Z}_q)\ncong 0,
\check{H}^2(S^3/X, S^3 - X; \mathbb{Q}) \ncong 0.
\end{equation}

Consider the following exact sequence of the pair $(S^3/X, S^3 - X; \mathbb{Z}_p)$:
$$\check{H}^{2}(S^3 - X; \mathbb{Z}_p) \stackrel{\delta}\longrightarrow
\check{H}^3(S^3/X, S^3 - X; \mathbb{Z}_p) \longrightarrow
\check{H}^3(S^3/X; \mathbb{Z}_p)\longrightarrow \check{H}^{3}(S^3 - X;
\mathbb{Z}_p).$$ 
The groups $\check{H}^{2}(S^3 - X; \mathbb{Z}_p)$ and
$\check{H}^{3}(S^3 - X; \mathbb{Z}_p)$  are trivial by (\ref{H^2(S3-x)})
and (\ref{H^3(S3-x)}), respectively. Next, observe that $\check{H}^3(S^3/X;
\mathbb{Z}_p) \cong \mathbb{Z}_p$. Therefore
\begin{equation}\label{cmZ_p3}
\check{H}^3(S^3/X, S^3 - X; \mathbb{Z}_p) \cong \mathbb{Z}_p.
\end{equation}
Let us show  that $S^3/X$ is  $clc^{\infty}_{Z_p}$ at all points. Obviously,
the space $S^3/X$ is a $clc^{\infty}_{Z_p}$ space for all points except at
the
point $x = X/X$ since  $S^3 -  X$ is an open manifold. As
 mentioned before, the sets $\{U_i/X\}_{i \in \mathbb{N}}$  form a
neighborhood base  of the point $x.$ Consider the groups
$\check{H}^n(U_i/X, X/X; \mathbb{Z}_p).$ By the Excision Axiom it follows that
$\check{H}^n(U_i/X, X/X; \mathbb{Z}_p) \cong \check{H}^n(U_i, X; \mathbb{Z}_p).$

From the following commutative diagram with exact rows
$$\begin{array}{llllll}
0 \cong \check{H}^{0}(X; \mathbb{Z}_p)& \longrightarrow &
\check{H}^1(U_i, X; \mathbb{Z}_p)& \longrightarrow & \check{H}^1(U_i; \mathbb{Z}_p) &\cong \mathbb{Z}_p\\
\downarrow &&\downarrow \pi^i &&&\downarrow \times n_i\\
0 \cong\check{H}^{0}(X; \mathbb{Z}_p)& \longrightarrow &
\check{H}^1(U_{i+1}, X; \mathbb{Z}_p)& \longrightarrow & \check{H}^1(U_{i+1}; \mathbb{Z}_p) &\cong \mathbb{Z}_p\\
\end{array}$$
it follows that for  a large enough $i,$  the homomorphism
$$\check{H}^1(U_i, X; \mathbb{Z}_p)\stackrel{\pi^i}\longrightarrow
\check{H}^1(U_{i+1}, X; \mathbb{Z}_p)$$is trivial. Therefore
\begin{equation}\label{clcZ_p1}  \  S^3/X \ \hbox{is a} \ 1-clc_{\mathbb{Z}_p} \hbox{space}. 
\end{equation}
By the analogous diagram for the group of coefficients $\mathbb{Z}$ 
it is easy to
see that the homomorphism $\check{H}^1(U_i, X;
\mathbb{Z})\stackrel{\pi^i}\longrightarrow \check{H}^1(U_{i+1}, X; \mathbb{Z})$ is
a monomorphism of the group $\mathbb{Z}$.  Therefore
\begin{equation}\label{notclcZ_p1}
S^3/X \ \hbox{ is not}\ 1-clc_{\mathbb{Z}}.
\end{equation}

By the exact sequence
$$\check{H}^1(X; \mathbb{Z}_p) \longrightarrow \check{H}^2(U_i, X; \mathbb{Z}_p)\longrightarrow
\check{H}^2(U_i, \mathbb{Z}_p)$$ since $\check{H}^2(U_i, \mathbb{Z}_p) \cong 0$ and the \v{C}ech cohomology group
$\check{H}^1(X; \mathbb{Z}_p)$ is obviously isomorphic to the direct limit of the
sequence
$$\mathbb{Z}_p \stackrel{\times n_1}{\rightarrow} \mathbb{Z}_p \stackrel{\times n_2}
{\rightarrow} \mathbb{Z}_p \stackrel{\times n_3}{\rightarrow}
\cdots$$ it follows that $\check{H}^2(U_i, X; \mathbb{Z}_p)\cong 0.$ By the Excision Axiom 
it follows that $$\check{H}^2(U_i/X, X/X; \mathbb{Z}_p)\cong 0$$ and 
\begin{equation}\label{clcZ_p2}
S^3/X \ \hbox{is a} \ 2-clc_{\mathbb{Z}_p} \hbox{space}. 
\end{equation}

By the following exact sequence of the pair $(U_i, X)$
$$\check{H}^2(X; \mathbb{Z}_p) \longrightarrow \check{H}^3(U_i, X; \mathbb{Z}_p)\longrightarrow
\check{H}^3(U_i, \mathbb{Z}_p)$$ and since the space $X$ is $1-$dimensional
and $U_i$ is homotopy equivalent to the circle it follows that $\check{H}^3(U_i, X; \mathbb{Z}_p)\cong 0$ therefore
$\check{H}^3(U_i/X, X/X; \mathbb{Z}_p)\cong 0$ and
\begin{equation}\label{clcZ_p3}
\  S^3/X \ \hbox{is a} \ 3-clc_{\mathbb{Z}_p} \hbox{space}. 
\end{equation}

By the local connectedness of the space $S^3/X$, by (\ref{clcZ_p1}),
(\ref{clcZ_p2}),
(\ref{clcZ_p3})
and since $\dim S^3/X =3$ it
follows that $S^3/X$ is a $clc^{\infty}_{\mathbb{Z}_p}$ space and $S^3/X$ satisfies the condition (2) of Definition \ref{chm}.

It follows by
(\ref{cmZ_p1}),
(\ref{cmZ_p2}), and
(\ref{cmZ_p3}) 
that $S^3/X$ satisfies the condition (3) of Definition \ref{chm}, therefore $S^3/X$ is
a $3-cm_{\mathbb{Z}_p}$ and a $3-hm_{\mathbb{Z}_p}$.

However, the space $S^3/X$ is neither
a
$3-cm_{\mathbb{Z}_q}$
nor
a
$3-cm_{\mathbb{Q}}$
since $\check{H}^2(S^3/X, S^3 - X; \mathbb{Z}_q)\ncong 0$ and $\check{H}^2(S^3/X, S^3 - X; \mathbb{Q})
\ncong 0  $
by
(\ref{notclcZ_p1}),
and is not
a
$3-cm_{\mathbb{Z}}$ since it is not
a
$1-clc_{\mathbb{Z}}.$
This completes the proof.\qed

\section{Epilogue}

The spaces which we have constructed are not ANR's, so the following remains an open question:

\begin{qws}
Let $\mathcal{P}$ be any set of  prime numbers. 
Does there exist a 3-dimensional ANR $X_{\mathcal{P}}$ with the following properties:
\begin{enumerate}
\item
$X_{\mathcal{P}}$
is a 
${\mathbb{Z}_p}$-homology 
3-manifold 
  for every
$p\in \mathcal{ P}$; and
\item
$X_{\mathcal{P}}$
 is not a 
${\mathbb{Z}_q}$-homology 3-manifold  for any $q\not\in
\mathcal{ P}$?
\end{enumerate}
\end{qws}


\section{Acknowledgements}

This research was supported by the Slovenian Research Agency
grants No. BI-US/13-14-027, P1-0292, J1-4144, J1-5435, and J1-6721.
We thank Allen Hatcher for noticing an inaccuracy in the initial version of Example 2.3.\\

\end{document}